\documentclass[fleqn,usenatbib]{mnras}

\usepackage{newtxtext,newtxmath}
\usepackage[T1]{fontenc}

\DeclareRobustCommand{\VAN}[3]{#2}
\let\VANthebibliography\thebibliography
\def\thebibliography{\DeclareRobustCommand{\VAN}[3]{##3}\VANthebibliography}

\usepackage{graphicx}	% Including figure files
\usepackage{amsmath}	% Advanced maths commands
\usepackage{subcaption}
\usepackage[dvipsnames]{xcolor,colortbl}
\usepackage{ulem}

\newcommand{\calE}{\mathcal{E}}

\definecolor{Gray}{gray}{0.85}
\newcolumntype{a}{>{\columncolor{Gray}}c}
\newcolumntype{?}{!{\vrule width 0.5pt}}
\newcommand{\be}{\begin{equation}}
\newcommand{\ee}{\end{equation}}
\newcommand{\bea}{\begin{eqnarray}}
\newcommand{\eea}{\end{eqnarray}}
\newcommand{\non}{\nonumber \\}

\definecolor{seagreen}{rgb}{0.190, 0.525, 0.361}

\def\bse{\begin{subequations}}
\def\ese{\end{subequations}}
\def\eps{\epsilon}
\def\presub{\vspace{.5cm} \noindent}
\def\bsig{{\bar \sigma}}
\title[The Three-Body Problem]{Measurement of three-body chaotic absorptivity predicts chaotic outcome distribution}
\author[Manwadkar et al.]{Viraj Manwadkar,$^{1,2}$\thanks{E-mail: virajmanwadkar@gmail.com}
Alessandro A. Trani,$^{3,4,5}$\thanks{E-mail: aatrani@gmail.com}
Barak Kol$^{6}$\thanks{E-mail: barak.kol@mail.huji.ac.il} %}
\\
$^{1}$ Department of Physics, Stanford University, 382 Via Pueblo Mall, Stanford, CA 94305, USA\\
$^{2}$ Kavli Institute for Particle Astrophysics \& Cosmology, P. O. Box 2450, Stanford University, Stanford, CA 94305, USA\\
$^{3}$Niels Bohr International Academy, Niels Bohr Institute, Blegdamsvej 17, 2100 Copenhagen, Denmark\\
$^{4}$Research Center for the Early Universe, School of Science, The University of Tokyo, Tokyo 113-0033, Japan\\
$^{5}$Okinawa Institute of Science and Technology, 1919-1 Tancha, Onna-son, Okinawa 904-0495, Japan\\
$^{6}$ Racah Institute of Physics, Hebrew University, Jerusalem 9190401, Israel\\
}

\date{Accepted XXX. Received YYY; in original form ZZZ}

\pubyear{2021}

\begin{document}
\label{firstpage}
\pagerange{\pageref{firstpage}--\pageref{lastpage}}
\maketitle

% Abstract of the paper
\begin{abstract}
The flux-based statistical theory of the non-hierarchical three-body system predicts that the chaotic outcome distribution reduces to the chaotic emissivity function times a known function, the asymptotic flux. Here, we measure the chaotic emissivity function (or equivalently, the absorptivity) through simulations. More precisely, we follow millions of scattering events only up to the point when it can be decided whether the scattering is regular or chaotic. In this way, we measure a tri-variate absorptivity function. Using it, we determine the flux-based prediction for the chaotic outcome distribution over both binary binding energy and angular momentum, and we find good agreement with the measured distribution. This constitutes a detailed confirmation of the flux-based theory, and demonstrates a considerable reduction in computation to determine the chaotic outcome distribution.
\end{abstract}

% This is a simple template for authors to write new MNRAS papers.
% The abstract should briefly describe the aims, methods, and main results of the paper.
% It should be a single paragraph not more than 250 words (200 words for Letters).
% No references should appear in the abstract.

% Select between one and six entries from the list of approved keywords.
% Don't make up new ones.
\begin{keywords}
chaos, gravitation, celestial mechanics, planets and satellites: dynamical evolution and stability
\end{keywords}
%%%%%%%%%%%%%%%%%%%%%%%%%%%%%%%%%%%%%%%%%%%%%%%%%%

%%%%%%%%%%%%%%%%% BODY OF PAPER %%%%%%%%%%%%%%%%%%

\section{Introduction}
\label{sec:intro}

The Newtonian three-body problem is one of the richest, most-fruitful and longest-standing open problems in physics. It is the fertile soil that grew numerous scientific theories including perturbation theory, the symplectic formulation of mechanics (Poisson brackets), (manifold) topology and chaos.   

Since Poincar\'e, a deterministic general solution is believed impossible \citep{Poincare1890}. Special cases that allow for analytic treatment are known to include the lunar limit, which displays orbit hierarchy, and the planetary limit that consists of a dominant mass, and hence a mass hierarchy. However, the general, non-hierarchical, case is believed to be non-amenable to an analytic deterministic treatment. Beginning with \citet{AgekyanAnosova1967}, computers were used to simulate this system and measure its outcome distribution. A statistical theory was formulated by \citet{monaghan76a}, and more recently, important advances were introduced to it by \cite{stone19} and \cite{ginat20}. However, this statistical theory suffers from certain shortcomings, including the introduction of a spurious parameter, the strong interaction region, which acts as a cutoff and a rough separator of regular and chaotic motion.

Inspired by \cite{stone19}, the flux-based statistical theory introduced in \cite{Kol2020} remedies this issue, and identifies probabilities not with phase-space volume, as in previous theories, but rather with phase-space flux. Its central result is the following reduction of the chaotic outcome distribution \be
 d\Gamma(u) = \frac{1}{\bsig}\, \calE(u)\, dF(u)
\label{dGamma}
\ee
where $u$ is a collective notation for outcome parameters; $d\Gamma(u)$ is the chaotic differential decay rate, namely the probability per unit time to decay from a chaotic state into an element $du$ of outcome parameters, and hence $d\Gamma(u)$ is proportional to the chaotic outcome distribution; $\bsig$ is the regularized chaotic phase-volume determined in \citet{sigma} ; $\calE(u)$ is the chaotic emissivity function, which is the probability of an outgoing state with outcome parameters $u$ to have originated from chaotic motion; finally, $dF(u)$ is the distribution of asymptotic flux. This relation constitutes a reduction since $dF(u)$ was determined in closed form, see \eqref{dF} below, and $\bsig$ is only a $u$-independent normalization, thereby reducing the study of  $d\Gamma(u)$ to that of $\calE(u)$.

Parts of the flux-based statistical theory were validated through simulations in \citet[]{manwadkar20} and \citet[]{Manwadkar21}, hereafter MKTL21. In addition, this theory stimulated a novel formulation of the three-body system, one that provides a natural dynamical reduction \citep{DynRed}.

\begin{figure*}
    \centering
    \includegraphics[width = \textwidth]{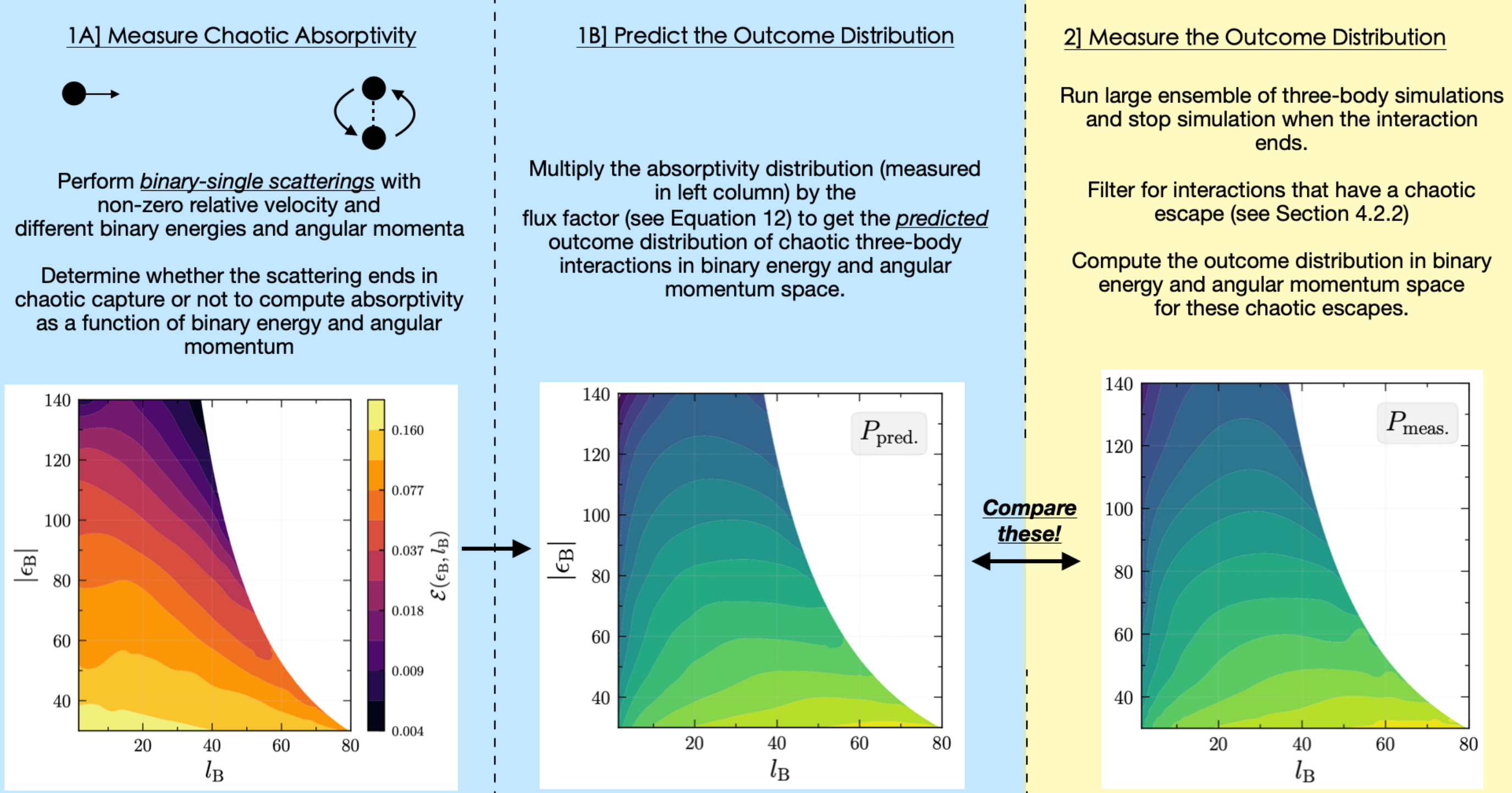}
    \caption{ The figure shows the general outline
    of this work. (\textit{Left Panel}) Description of our method to measure chaotic absorptivity from binary-single scattering interactions in different configurations. (\textit{Center Panel}) Description of procedure to predict the outcome distribution of chaotic three-body interactions that uses the chaotic absorptivity measurement from left panel. (\textit{Right Panel}) Description of procedure to independently measure the chaotic three-body outcome distribution by means of three-body interaction simulations. The goal of the paper is to make a comparison between this measured outcome distribution and the predicted outcome distribution to test the validity of the flux-based theory presented in \citet{Kol2020}.     }
    \label{fig:pap_schema}
\end{figure*}

The goal of this paper is to test the main reduction \eqref{dGamma} of the flux-based theory, by measuring the chaotic emissivity function $\calE(u)$. Through time reversal symmetry, $\calE$ can be interpreted also as the absorption function, namely the probability of an ingoing state with income parameters $u$ to proceed to chaotic motion. In other words, in analogy with Kirchhoff's law of thermal radiation, chaotic absorptivity and emissivity are identical. Once the chaotic absorptivity function is measured, the chaotic outcome distribution is predicted and compared with the its measured distribution. The road taken in the paper is shown schematically in Fig. \ref{fig:pap_schema}.

A measurement of chaotic absorptivity requires to simulate evolutions of incoming states only until the three bodies meet (and shortly thereafter), while the measurement of outcome distribution requires to evolve the system until it disintegrates, which is typically much longer. In this sense, \eqref{dGamma} provides a meaningful computational reduction of the measurement of the chaotic outcome distribution. Future analytic expressions or approximations for $\calE$ may provide further reduction.

This paper is organized as follows. We start in Section \ref{sec:theory} by setting up the problem and by reviewing the flux-based theory. Section \ref{sec:method} describes our method of measurement and Section \ref{sec:sims} describes the simulations. Data analysis an results are presented in Section \ref{sec:data}. We conclude with a summary and discussion in Section \ref{sec:summary}.

\section{Statistical Theory}
\label{sec:theory}
% BK: A review of the flux-based reduction of the outcome statistics

In this section, we briefly setup the problem and review the flux-based theory introduced in \citet{Kol2020}. Then we recall and discuss the predictions to be tested in the current paper.

\presub {\bf Setup}. The three-body problem can be defined though the Hamiltonian \be
H\left( \{\vec{r}_c,\, \vec{p}_c \}_{c=1}^3 \right) :=  \sum_{c=1}^3 \frac{\vec{p}_c^{~2}}{2 m_c} - \left( \frac{G\, m_1\, m_2}{r_{12}} + cyc. \right)  \label{def:H} 
\ee
where $m_c, ~c=1,2,3$ denote the three masses, $\vec{r}_c$  the bodies' position vectors and $\vec{p}_c$  their momenta, $G$ is Newton's gravitational constant, and $r_{cd}=\left| \vec{r}_c - \vec{r}_d \right|$. 

The conserved charges are the total linear momentum, the total energy, and the total angular momentum denoted by $\vec{P},\, E,\, \vec{L}$, respectively. We work in the center of mass frame.

The conservation laws allow a disintegration of the system. Generally, a non-hierarchical three-body motion ends in this way. Assuming negative total energy, the components of the outgoing states are a binary and a single, also known as the escaper. When they are far apart, the system decouples into a binary subsystem, and an effective hierarchical system defined by replacing the binary by a fictitious object obtained by collapsing it to its center of mass. The effective system describes the relative motion of the binary and the single. The outcome parameters $u$ are given by \be
 u=\left( s,\, \eps_B,\, \vec{l}_B,\, \psi_B,\, \eps_F,\, \vec{l}_F,\, \psi_F \right)
\ee
where $s=1,2,3$ is the escaper identity, $\eps_B,\, \vec{l}_B,\psi_B$ denote the energy, angular momentum and pericenter angle (measured relative to the line of nodes) of the decoupled binary, and $\eps_F,\, \vec{l}_F,\psi_F$ denote the analogous quantities for the effective system. These variables obey the obvious relations $\eps_B+\eps_F=E,~ \vec{l}_B + \vec{l}_F=\vec{L}$. Hence there are altogether 6 independent continuous outcome parameters. Symmetry with respect to rotations around $\vec{L}$, means that $u$-dependent quantities depend essentially only on 5 outcome parameters.
% include the energy and angular momentum of the decoupled binary, denoted by $\eps_B,\, \vec{l}_B$, and those of the effective system, denoted by $\eps_F,\, \vec{l}_F$. 

A special role is played by the binary constant defined by \be
 k := \frac{m_a\, m_b}{m_a + m_b} \left( G m_a m_b\right)^2 \label{def:k}
 \ee
 where $m_a,\, m_b$ are the masses which compose the binary. $k_s$ denotes the binary constant of the binary defined by an escaper $s=1,2$ or $3$.

For a more complete setup,  see \citet{Kol2020}. 

\presub {\bf Flux-based statistical theory}. The flux-based theory differs from all past statistical treatments starting with \citet{monaghan76a} and up to the review book \citet{Valtonen_book_2006} and the closed-form determination of the outcome distribution \cite{stone19}. All previous treatments assume the micro-canonical ensemble, namely, assign probabilities according to phase-space volume. Moreover, they introduce the so-called strong interaction radius in order to guarantee finite phase-space volumes and to exclude an irrelevant part of phase-space (which correspond to causally inaccessible escape scenarios). The use of the micro-canonical ensemble implicitly assumes a closed system with a bounded phase space, and random (or technically, ergodic) motion, whereas the three-body system has an unbounded phase space, is open to disintegration, and its phase space is divided between regular and chaotic motion. In addition, setting the value of the strong interaction radius is somewhat arbitrary. 

% Figure of binary + effective systems ?

The flux-based theory \citep{Kol2020} remedies these flaws by focusing on phase-volume flux, rather than the phase-volume itself. This is natural for a disintegrating system. The flux is inherently finite and independent of the location where it is measured, and in this way infinite probabilities never appear, and a cutoff is no longer necessary.  

We note that the time evolution of a three-body system can be divided into three kinds of motion: asymptotic motion where at least one of the bodies is far away, chaotic interaction, and finally sub-escape excursions, which are a part of motion where the system clearly separates into a binary and a single which fly away from each other, yet their relative velocity is below the escape velocity, and hence they are bound to fall back towards the center of mass, and interact. Accordingly, the sub-escape excursions can be called quasi-asymptotic states. This decomposition of phase space is similar, if not identical, to ones found in the literature. In particular, this paper's ``regular scattering'' is similar to ``non-resonant interaction'' in the language of \citet{McMillanHut1996}, ``chaotic interaction'' resembles   ``democratic resonance'' there, and sub-escape excursion resembles ``hierarchical resonance'' there.
 
The first part of the flux-based theory considers the system's probability distribution as a time-dependent variable. The probability measure is divided between the ergodic region, the sub-escape excursions and the asymptotic states, and the latter two are further continuously distributed over their parameters. This distribution differs from the micro-canonical ensemble of previous treatments. The time evolution of the distribution is formulated through the system of equations (2.35) of \citet{Kol2020}, whose solution describes the statistics of outcome parameters and decay times. 

The second part of the theory involves the differential decay rate out of the ergodic region, which is an essential ingredient in the equation for the above-mentioned statistical evolution of the system. Its distribution over all the asymptotic state parameters $u$ is denoted by $d\Gamma(u)$, and it can be shown to factorize exactly according to \eqref{dGamma} where the distribution of asymptotic flux, $dF(u)$ is given by \be
	dF(u) = 2\pi\, \frac{\sqrt{k_s}\, d\eps_B}{(-2 \eps_B)^{3/2}} \, \frac{1}{l_B\, l_F} d^3l_B\, d^3 l_F\, \delta^{(3)}\left(\vec{l_B}+\vec{l_F}-\vec{L}\right)\, d\psi_B\, d\psi_F
%	d\Gamma_s(u) \propto  \calE (u) \cdot \frac{\sqrt{k_s}\, d\eps_B}{(-2 \eps_B)^{3/2}} \, \frac{1}{l_B\, l_F} d^3l_B\, d^3 l_F\, \delta^{(3)}\left(\vec{l_B}+\vec{l_F}-\vec{L}\right)
\label{dF}
\ee
where $s=1,2,3$  is the escaper identity, $\calE$ is the chaotic emissivity ($=$ absorptivity) and $k_s$ is defined in (\ref{def:k}). The variables $\eps_B, l_B$ range over the domain \be
-2\, \eps_B\, l_B^2  \le  k_s ~, \qquad
 \eps_B  \le  E  ~.
\label{binary_region}
\ee

Chaotic absorptivity $\calE(u)$ is defined to be the probability that a scattering of a single off a binary with a random mean anomaly (time along orbit) will evolve into a chaotic trajectory, rather than a regular scattering such as a flyby or a regular exchange, which leads to prompt ejection.
%, and the scattering parameters are specified by $u$, see eq. 2.23 of \citet{Kol2020}. 
%In analogy with Kirchhoff's law of thermal radiation, the chaotic absorptivity with incoming parameters $u$ equals chaotic emissivity with outgoing parameters $u$. 
Clearly, $\calE(u)$ serves to account for the division of phase space into regions of regular and chaotic motion. We shall often omit the adjective ``chaotic'' since we do not discuss other types of emissivity.   

%The rest of the expression accounts for the flux of phase space volume into the asymptotic states. This factor 
The asymptotic flux $dF(u)$
reflects the fact that this theory is based on the framework of an open, rather than closed, chaotic system, where the flux of phase-space volume replaces the volume itself as a measure of probability. Hence, it is called the flux-based theory. To have some insight into the form of the flux factor, note that $2 \pi \sqrt{k_s}/(-2 \eps_B)^{3/2}$ is the binary period, and the denominator factors $l_B,\, l_F$ each originate in the central force nature of the respective systems. 

By definition, $\calE(u)$ is bounded to the range $0 \le \calE (u) \le 1$. Otherwise, so far, it is an unknown function of the asymptotic parameters $u$. In this way, (\ref{dGamma}) factors out the outgoing flux and reduces the determination of $d\Gamma(u)$ to that of $\calE(u)$. For derivations, see \citet{Kol2020}.

% Altogether, through these ingredients, the flux-based theory addresses and improves upon the above-mentioned issues of the previous statistical theories. 

\section{Measurement Method}
\label{sec:method}

% Marginalization
% Criterion for absorption
% Parameter values

In this section, we describe several aspects of the measurement method: marginalization, the criterion for absorption, and the parameter values that we chose.

\presub {\bf Marginalization}. The differential decay rate \eqref{dGamma} is distributed over the 6d space of outcome parameters $u$. In practice, we measure its distribution over a smaller set of parameters $v$, so that \be
 u=(v,w)
 \ee
 where $w$ denotes the remaining outcome parameters.  
 
The marginalized decay rate is defined by \be
\frac{d\Gamma}{dv}(v) := \int \frac{d\Gamma}{dv dw}(v,w)\, dw
\ee
The reduction \eqref{dGamma} implies that it can be written as \bea
    \frac{d\Gamma}{dv}(v) &=& \frac{1}{\bsig} \int \calE(v,w)\, \frac{dF}{dv dw}(v,w)\, dw \non
            &=& \frac{1}{\bsig}\, \calE(v)\, \frac{dF}{dv}(v)
\eea
where the marginalized asymptotic flux and absorptivity are defined by \bse \begin{align}
    \frac{dF}{dv}(v) &:= \int \frac{dF}{dv dw}(v,w)\, dw \label{marg_flux_gen} \\
    \calE(v) &:= \frac{1}{dF/dv} \int \calE(v,w)\, \frac{dF}{dv dw}(v,w)\, dw \label{marg_calE_gen} 
    \end{align}
\ese
%\bea
%    \frac{dF}{dv}(v) &:=& \int \frac{dF}{dv dw}(v,w)\, dw \non
%    \calE(v) &:=& \frac{1}{dF/dv} \int \calE(v,w)\, \frac{dF}{dv dw}(v,w)\, dw
%\eea

The last relation will be used in the next section  to define the marginalized, or averaged absorptivity. In words, it means that {\it this averaging of $\calE$ is weighted according to asymptotic flux}. Note that in addition to the above mentioned marginalization, by definition, a measurement of $\calE$ requires an average over binary phase (more precisely, over the mean anomaly) \citep{Kol2020}.

In practice, we consider two specific marginalizations. The first marginalization is is over the pericenter angles $w=(\psi_B,\psi_F)$. In this case, the flux weighting amounts to $\int d\psi_B/(2\pi) \, d\psi_F/(2\pi)$, namely, uniform weight. The remaining outcome parameters are $v=(\eps_B,\, \vec{l}_B)$. 

In the second marginalization, we further marginalize over the angle $\sigma$ between $\vec{l}_B$ and $\vec{L}$, so that the $\calE$ is distributed only over the magnitude $l_B$, namely $v=(\eps_B,\, l_B)$. The integration measure over $\sigma$ is given by \be
 \frac{d\cos \sigma}{l_F} \equiv - \frac{dl_F}{L\, l_B}
\label{shell_differential}
\ee
where $l_F^2(\sigma)=L^2+l_B^2-2\, L\, l_B\, \cos \sigma$. In particular, we record for later use an explicit expression for the general marginalized flux \eqref{marg_flux_gen} \be
 \frac{dF}{d\eps_B\, dl_B} \propto \frac{l_B}{(-\eps_B)^{3/2}} ~.
\label{marg_flux}
\ee
This is gotten as follows $dF \propto d\eps_B/(-\eps_B)^{3/2}\,  d^3l_B/(l_B\, l_F) = d\eps_B/(-\eps_B)^{3/2}\, l_B\, dl_B \int d^2\Omega/l_F = d\eps_B/(-\eps_B)^{3/2}\, l_B\, dl_B 4 \pi/L \propto d\eps_B\, l_B\, dl_B/(-\eps_B)^{3/2}.$ The second to last expression is gotten by Newton's shell theorem (related to the identity \eqref{shell_differential}), which relies on the inequality $l_B \le L$, valid for the parameter values of the case at hand (detailed later in this section).

\presub {\bf Criterion for chaotic absorption}. This criterion in described in full in Section \ref{sec:sims} (in the second paragraph of Section~\ref{sssec:sims_absorb} and below Eq. \ref{rbarak}). Here we include motivation and an informal description.

On general grounds, the mapping of incoming states to outgoing states defined by a three-body system contains both regularity islands and chaotic regions. A study of the phase portraits of the three-body system, such as Fig.~3 of \citet{manwadkar20} revealed to us two kinds of islands: one associated with a single close encounter, and another associated with a sequence of two close encounters.

We define a close encounter as a state where two of the bodies are near each other, so there exists a hierarchy of distances, while at the same time, the energy of the two in their center of mass frame is positive, namely, the two are unbound. 

This observation suggests to identify a time evolution as a chaotic absorption if it is followed by 3 or more close encounters before one of the bodies is ejected. However, we found the number of sufficiently equilateral, or democratic, configurations, to be a more robust observable than the number of close encounters. Generally, the number of democratic configurations before ejection is the successor of the number of close encounters (since every time evolution starts with a democratic configuration, and after it, close encounters and democratic configurations generally alternate). Therefore, we have identified a time evolution as a chaotic absorption if it results in 4 or more democratic configurations before one of the bodies is ejected into either an escape or an excursion. This can be summarized informally by \bea
 \mbox{no. of democratic config's before 1st ejection} \ge 4 \non
 \implies \mbox{chaotic absorption}
\eea

\presub {\bf Parameter values}. We choose the mass parameters to be equal, which is the most symmetric choice, and hence a good place to start at. The masses are set to $m_1=m_2=m_3=15$ in $\mathrm{M}_\odot$ units, or equivalently, in $N$-body mass units.

The conserved charges correspond to a class of states composed of a circular binary and a far away tertiary at rest. These are the initial conditions used in \citet{manwadkar20,Manwadkar21} to measure the outcome distribution. In $N$-body units, where Newton's constant is set to unity $G=1$, we have $E=-27, L=75 \sqrt{3/2} \simeq 91.86$. This correspond to a circular binary with radius $5$ simulation length units, whose center of mass is at a distance of approximately $100$ from the tertiary. The binary constants, $k_s$ \eqref{def:k}, are $s$-independent here, namely $k_s=k$, and the conserved charges define a dimensionless parameter $2|E|L^2/k=1.2$.

\section{Simulations}
\label{sec:sims}

\subsection{TSUNAMI $N$-Body code}
\label{ssec:tsunami}

As in MKTL21, we run the three-body simulations with the regularized $N$-body code \textsc{tsunami} \citep{trani2019,trani2022b}. The main advantage of \textsc{tsunami} over other integrators is that it implements the algorithmic regularization scheme of \citet{mikkola99a,mikkola99b}. This scheme increases the dramatically the accuracy of the integration, making sure that the simulations do not stall or lose accuracy during close encounters. Together with the Bulirsch-Stoer extrapolation scheme \citep{stoer80} and the chain-coordinate system \citep{mikkola1993}, it makes \textsc{tsunami} especially suited to model few-body interactions. For more details on the code, we refer to MKTL21 and \citep{trani2022b}. Because here we focus only on Newtonian dynamics of point-masses, we disable additional forces like % BK omit: the 1PN, 2PN and 2.5PN 
post-Newtonian corrections \citep{blanchet14} and tides \citep{hut81,samsing18_tides}. Likewise, we neglect collisions between particles.

To determine the state of the three-body system, we adopt a similar classification scheme as that employed in \citet{manwadkar20} and MKTL21. The hierarchy state of the triple is checked at every timestep by selecting the most bound pair and checking its relative energy with respect to the third body. If the third body is bound to the binary, the triple might be undergoing an excursion, during which the binary is relatively unperturbed by the single which was ejected with a speed below the escape velocity. To determine this, we first check whether the binary is perturbed by the single. We do this by comparing the relative force of the binary and the tidal force from the single using the following dimensionless quantity:
\begin{equation}\label{eq:ffactor}
	f_{\rm tid} = \frac{F_{\rm tid}}{F_{\rm rel}} = \frac{2 m_{\rm bin} m_3}{m_1 m_2} \left(\frac{a_{\rm bin} (1+e_{\rm bin})}{R}\right)^3 \,,
\end{equation}
where $m_{\rm bin} = m_1 + m_2$ is the total mass of the binary, $m_3$ is the mass of the single, $R$ is the distance of the single from the binary center-of-mass, and $a_{\rm bin}$, $e_{\rm bin}$ are semimajor axis and eccentricity of the binary, respectively. When $f_{\rm tid} \geq 1$, the tidal force is greater than the relative force of the binary at apocenter, and we deem the interaction as non-hierarchical. Conversely, if $f_{\rm tid} < 1$, we label the interaction as a candidate hierarchical excursion, which concludes when $f_{\rm tid} \geq 1$ again. We then compare the excursion duration $t_{\rm ex}$ with the median binary period $\langle P_{\rm in} \rangle$ during the candidate excursion, and label the interaction as an excursion only if the former is longer than the latter ($t_{\rm ex} > \langle P_{\rm in} \rangle$).

We also monitor the hierarchy of the system by using a modified version of the homology radius $R_{\cal H}$ used in \citet{manwadkar20}:
\begin{equation}\label{rbarak}
	\tilde{R}_{\cal H} = \dfrac{3 r_{\rm min}^2}{\sum_{i\neq j} r_{i,j}^2}
\end{equation}
for $i,j = 1,2,3$, and where $r_{i,j}$ is the relative distance between particle $i$ and $j$, and $r_{\rm min}$ is the minimum interparticle distance, i.e. $r_{\rm min} = \min{\{r_{i,j}\}}$.
This quantity runs in the range $0 \le \tilde{R}_{\cal H} \le 1$: the limit $\tilde{R}_{\cal H}=1$ corresponds to equidistant configurations (equilateral), while $\tilde{R}_{\cal H}=0$ corresponds to hierarchical configurations. 
% is close to one when the particles are roughly equidistant, i.e. when they are undergoing a democratic resonance, and close to zero when the system is hierarchically well separated. 
The modified homology radius $\tilde{R}_{\cal H}$ has the advantage of being a more smooth function of the triangle geometry as compared with $R_{\cal H}$, since it does not contain the ${\rm max}$ function in the denominator (this means that $\tilde{R}_{\cal H}$ is smooth at hierarchical isosceles configurations.)

Each time the value of $\tilde{R}_{\cal H}$ grows above a threshold chosen as $\tilde{R}_{{\cal H}, {\rm D}} = 0.33$ and goes back below $\tilde{R}_{{\cal H}, {\rm D}}$ again is counted as a
sufficiently equilateral (or Lagrangian) configuration, or in short a democratic configuration. 
%democratic resonance, or scramble. 
We denote by $N_{\rm D}$ the total number of democratic configurations occurring over an interaction.

\subsection{Initial setup and decision schemes}

In order to compare the measurement of the chaotic absorptivity with the differential decay rate, we perform two different kinds of simulations, using different initial setup and stopping conditions. 

\subsubsection{Measurement of chaotic absorptivity}
\label{sssec:sims_absorb}

\begin{figure}
    \centering
    \includegraphics[width = \columnwidth]{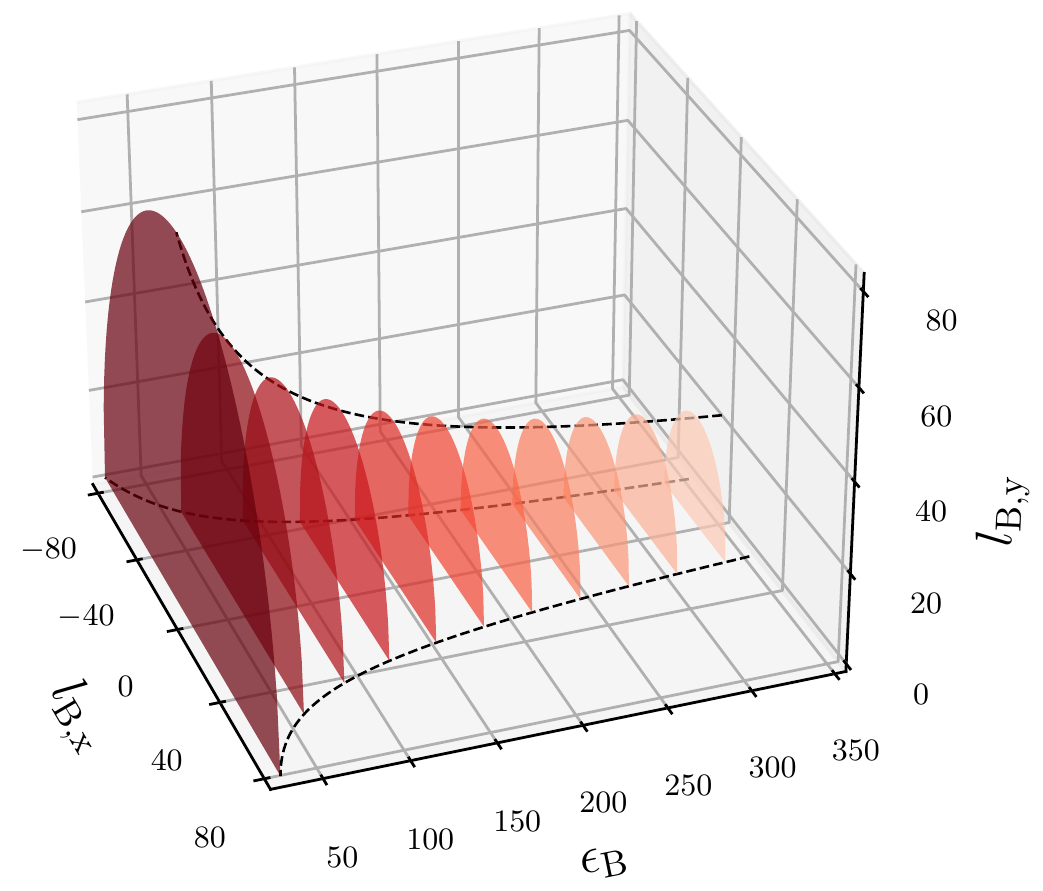}
    \caption{A 3D schematic representation of the outcome distribution of three-body interactions. Each semi-circular slice in this figure corresponds to a slice with constant binary energy $|\epsilon_{\rm B}|$. These slices are equivalent to the slices in Figure~\ref{fig:measured_calE}.}
    \label{fig:3dsec}
\end{figure}

For the chaotic absorptivity, the initial configuration is a hyperbolic encounter between a binary and a single body. For every set of simulations, we keep the total energy $E$ and angular momentum $L$ constant to the values mentioned at the last part of Section~\ref{sec:method}, and we probe the parameter space of binary energy $\epsilon_B$ and angular momentum $l_B$ in a grid-like fashion. We sample 100 values of $\epsilon_B$ on a regular spacing between $-30$ and $-150$, and select $l_B = 1.5, 2.5, 7.5, 10, 20, 30, 40, 50, 60, 70$. Note that the $\epsilon_B$--$l_B$ space has a forbidden region whose boundary corresponds to circular orbits (see Figure~\ref{fig:3dsec}) , for which
\begin{equation}
l_{B,\mathrm{max}}(\epsilon_B) = \sqrt{\frac{-k}{2 \epsilon_B}}
\label{eqn:max_lb_eqn}
\end{equation}
(see also Eq. \ref{binary_region}). We sample this boundary region by including sets of simulations with $l_{B,\mathrm{max}}$ for each of the $\epsilon_B$ considered. Our choice of initial conditions is summarized in Table~\ref{tab:ic1}. 
In our chosen reference frame, the binary lies on the $x$-$y$ plane, with the pericenter along the $x$ axis. Therefore, we set the argument of pericenter $\psi_B$, the inclination $i_B$ of the binary orbit to zero.
We randomize the mean anomaly of binary $M_b$, the argument of the pericenter of the hyperbolic orbit $\psi_F$, and the angular momentum $\vec{l}_F$ of the hyperbolic orbit. Specifically, the magnitude $l_F$ is drawn uniformly in the interval $[L - l_B, L + l_B]$ according to \eqref{shell_differential}, while the inclination of the hyperbolic orbit is determined so to conserve $l_B$ and $L$. The longitude of the ascending node of the hyperbolic orbit is then drawn uniformly in $[0, 2\pi]$. The initial position of the single with respect to the center of mass of the binary is set to 20 times the initial binary semi-major axis.

To determine whether an interaction was chaotically absorbed or not, we run the simulations until an excursion occurs. We then count the number of democratic configurations $N_{\rm D}$, according to our criterion described above. If $N_{\rm D} \geq 4$, we consider the interaction absorbed.

\begin{figure}
    \centering
    \includegraphics[width = \columnwidth]{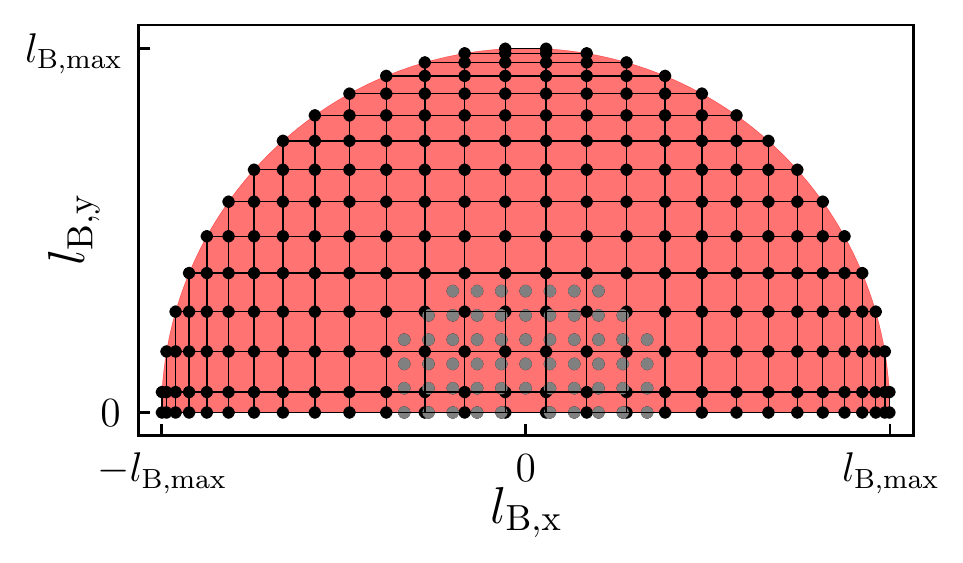}
    \caption{The grid in $l_{\rm B,x} - l_{\rm B,y}$ space (at a fixed $\epsilon_{\rm B}$) we use for measuring chaotic absorptivity. At each of these grid points, we run an ensemble of binary-single scattering experiments to measure the chaotic absorptivity at that given grid point. We use two super-imposed grids for better interpolation across grid. The Chebyshev grid is denoted by black circles (and with grid lines for reference). The smaller, uniform grid is denoted by grey circles. Refer to Section~\ref{sssec:sims_absorb} for details on the grid. }
    \label{fig:grid_layout}
\end{figure}

The sets described above are designed to measure the bi-variate chaotic absorptivity in the $l_{\rm B}$-$\epsilon_{\rm B}$ space. In addition, we also run a grid of simulation sets by fixing also the angular momentum of the hyperbolic orbit $l_F$. This can be also interpreted as fixing the angle between $\vec{l}_B$ and $\vec{l}_F$ (see Section~\ref{ssec:bivariate}). In this way, we can perform a measurement of the tri-variate absorptivity in terms of $l_{\rm B,x},\, l_{\rm B,y}$, and $\epsilon_{\rm B}$. For this set, we select $\epsilon_B$ in between $-30$ and $-160$ in equal spacing of $10$, and from $-160$ to $-300$ in equal spacing of $20$. For each value $\epsilon_B$ we sample a grid in $l_{B,x}$--$l_{B,y}$ space, where $l_{B,x}$ is the component of $\vec{l}_B$ aligned with $\vec{L}$. Using the angle between $\vec{l}_B$ $\vec{L}$, which we call $\sigma$  (see Section~\ref{sec:method}), the components of $l_B$ can be written as $l_{B,x} = l_B \cos\sigma$ and $l_{B,y} = l_B \sin\sigma$. The grid in $l_{B,x}$--$l_{B,y}$ is constructed using Chebyshev nodes with a super-imposed uniform grid of smaller radius (see Figure~\ref{fig:grid_layout}). The motivation for using two super-imposed grids is to get a better accuracy in grid interpolation at both the center and edges of disk. This choice of grid was found to be efficient for interpolation on a disk while considering the number of grid points needed and the corresponding accuracy. Inspired by the 1D Chebyshev grid, we define the 2D Chebyshev-like grid by first considering points on a semi-circle $(x_i,y_i)$ defined below and then constructing the grid by connecting all these points in a grid-like fashion as illustrated in Figure~\ref{fig:grid_layout}. The points on the semi-circle $(x_i,y_i)$ are defined as
\begin{align}
    \left( x_i,y_i \right) \in \Biggl\{ \left( l_{\rm B,max}\cdot \cos\left(\frac{2i -1}{2N} \pi\right) ,l_{\rm B,max} \cdot \sin\left(\frac{2i -1}{2N} \pi\right) \right) \Biggr\} \\
    + \Biggl\{(-l_{\rm B,max},0 ), (l_{\rm B,max},0) \Biggr\} 
\end{align}
for $i \in \{ 1,2,\cdots N \}$ for an even integer $N$. The smaller uniform grid is over a disk of radius $r = 0.4 \cdot l_{\rm B,max}$. We choose a value of $N=28$ for our 2D Chebyshev grid.

% \BK{Viraj, did you mean the following? \\
% Inspired by the 1d Chebyshev grid, we define the 2d Chebyshev-like grid by 
% \begin{equation}
%     \left( x_i,y_j \right) = \left( l_{\rm B,max}\cdot \cos\left(\frac{2i -1}{2N} \pi\right) ,l_{\rm B,max} \cdot \sin\left(\frac{2j -1}{2N} \pi\right) \right)
% \end{equation}
% for $i \in \{ 1,2,\cdots,N \}, j \in \{ 1,2,\cdots,N/2 \}$}

\subsubsection{Measurement of outcome distribution}
\label{sssec:outcome_meas}

For the outcome distribution, we adopt the same initial conditions as MKTL21. Namely, we simulate binary--single encounters where the single is initially at rest with respect to the binary center of mass (hence an impact parameter is undefined). The only two free parameters that we allow in this setup are the orbital phase of the binary and the inclination of the binary with respect to the line joining the binary center of mass and the single. The orbital phase is uniformly sampled between 0 and $2\pi$, while the inclination is uniformly sampled in the cosine between $-1$ and $1$. We simulate a total of $10^7$ realizations. Unlike in the numerical experiments for the measurement of the chaotic absorption, we run the simulations until the final breakup occurs. Our initial conditions for this set of simulation is summarized in Table~\ref{tab:ic2}. 

\begin{figure}
    \centering
    \includegraphics[width = \columnwidth]{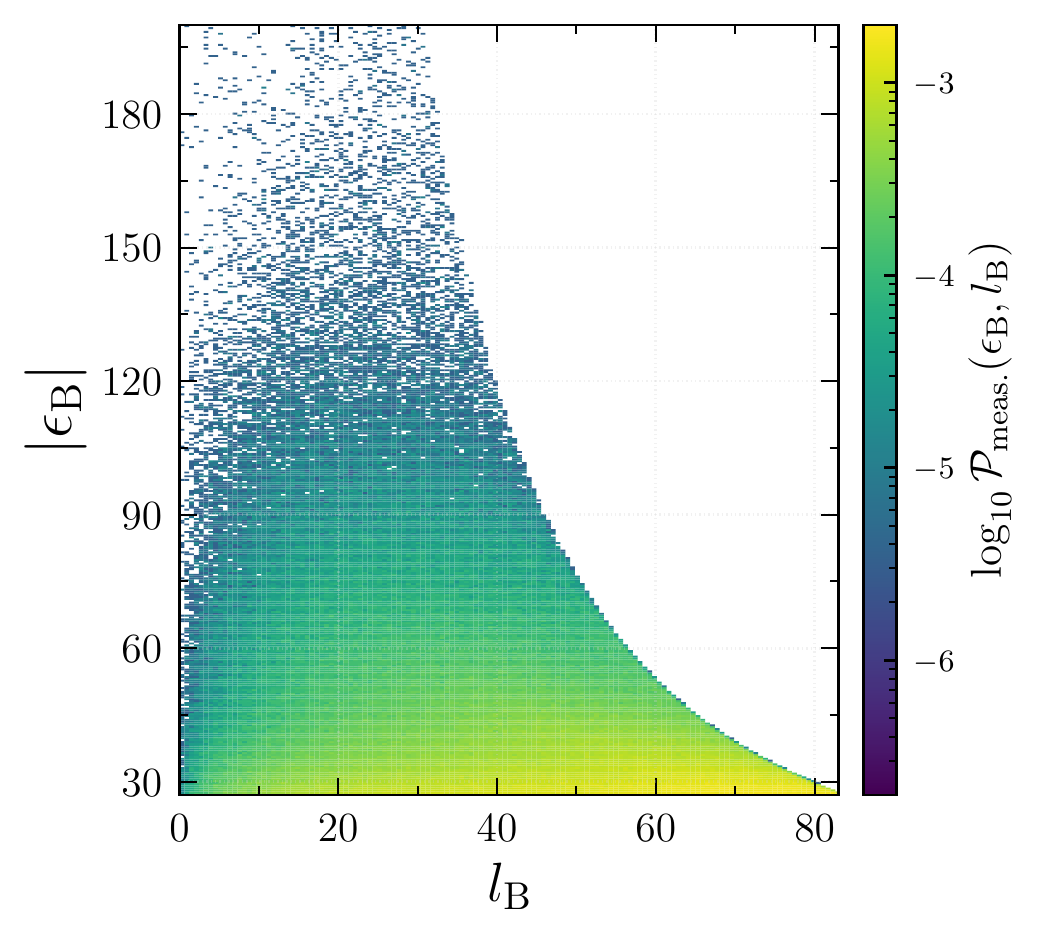}
    \caption{The 2D histogram distribution of interactions with a chaotic escape in $\eps_{\rm B} - l_{\rm B}$ space. Refer to Section~\ref{sssec:outcome_meas} for details on the cuts used to consider interactions with chaotic escapes. }
    \label{fig:outcome_hist}
\end{figure}

As we concern ourselves with the chaotic outcome distribution in this paper, we have to apply appropriate cuts to this ensemble of $10^7$ realizations to isolate the chaotic escapes. In line with analysis done in \citet{manwadkar20} and MKTL21, we first apply a cut on the lifetime of the total three-body interaction to remove the short-lived, regular interactions. Such cuts help in removing the regular islands in phase space (see Figure 3 of \citealt{manwadkar20}). Specifically in this work, we apply a cut of $\tau_{\rm lifetime} > 50$ yrs. Furthermore, to ensure a chaotic escape, we apply a cut to only consider interactions which have $\geq 4$ number of democratic configurations right before ejection (but after last excursion). The definition of a democratic configurations is the same as the one used in the binary-single scattering experiments described in Section~\ref{ssec:tsunami}. With these cuts to isolate chaotic escapes, we are left with $\sim 11\%$ of the total $10^7$ realizations. Figure~\ref{fig:outcome_hist} shows the 2D histogram distribution of binary energy $\eps_{\rm B}$ and angular momentum $l_{\rm B}$ of these interactions with chaotic escapes. To check that our sample of chaotic escapes is reasonable, we look at the ejection likelihoods for the 3 equal masses. In the chaotic escape limit, the ejection probabilities for the 3 equal masses should be equal, that is, 1/3. With these cuts, we find ejection probabilities of 0.332, 0.334, and 0.333 for the 3 masses. 

\begin{table}
\centering
\begin{tabular}{l|c|c}
      & Bi-variate $\calE \left( l_{\rm B},\, \epsilon_{\rm B} \right)$ &  Tri-variate $\calE \left( l_{\rm B,x},\, l_{\rm B,y},\, \epsilon_{\rm B} \right)$ \\\hline
     $\epsilon_B$ & ${\cal U}[-30, -150]$ & Spaced in $[-30, -300]$ \\
    $l_B$ & Spaced in $[1.5, 70]$ & Sampled in $[0, l_{B,\mathrm{max}}] $\\ 
    $l_F$ & ${\cal U}(L - l_B, L + l_B)$ & Chebyshev grid in $l_{B,x}$--$l_{B,y}$ \\
\end{tabular}
\caption{Initial conditions of our sets of simulations for the measurement of the chaotic absorptivity. $\epsilon_B$: energy of the binary in $N$-body units. $l_B$: angular momentum of the binary in $N$-body units. $l_F$: angular momentum of the binary-single hyperbolic orbit in $N$-body units. ${\cal U}(a,b)$: random uniform distribution between $a$ and $b$. For each set we run a total  of $10^5$ realizations. See Section~\ref{sssec:sims_absorb} for more details.}
\label{tab:ic1}
\end{table}

\begin{table}
\centering
\begin{tabular}{c|c|c|c}
    $m$ [M$_\odot$]  & $a_{\rm bin}$ [au] & $\lambda$ [rad] & $\cos{i}$ \\\hline
     15 & 5 & ${\cal U}(0, 2\pi)$ & ${\cal U}(1, -1)$\\
\end{tabular}
\caption{Initial conditions of our set of simulations for the measurement of the outcome distribution. $m$: mass of the three bodies. $a_{\rm bin}$: semimajor axis of the binary. $\lambda$: true longitude of the binary. $i$: inclination of the binary with respect to the line joining the center of mass of the binary and the single. ${\cal U}$: random uniform distribution between $a$ and $b$. We run a total of $10^7$ realizations. See Section~\ref{sssec:outcome_meas} for more details.}
\label{tab:ic2}
\end{table}

\begin{figure*}
    \centering
    \includegraphics[width = \textwidth]{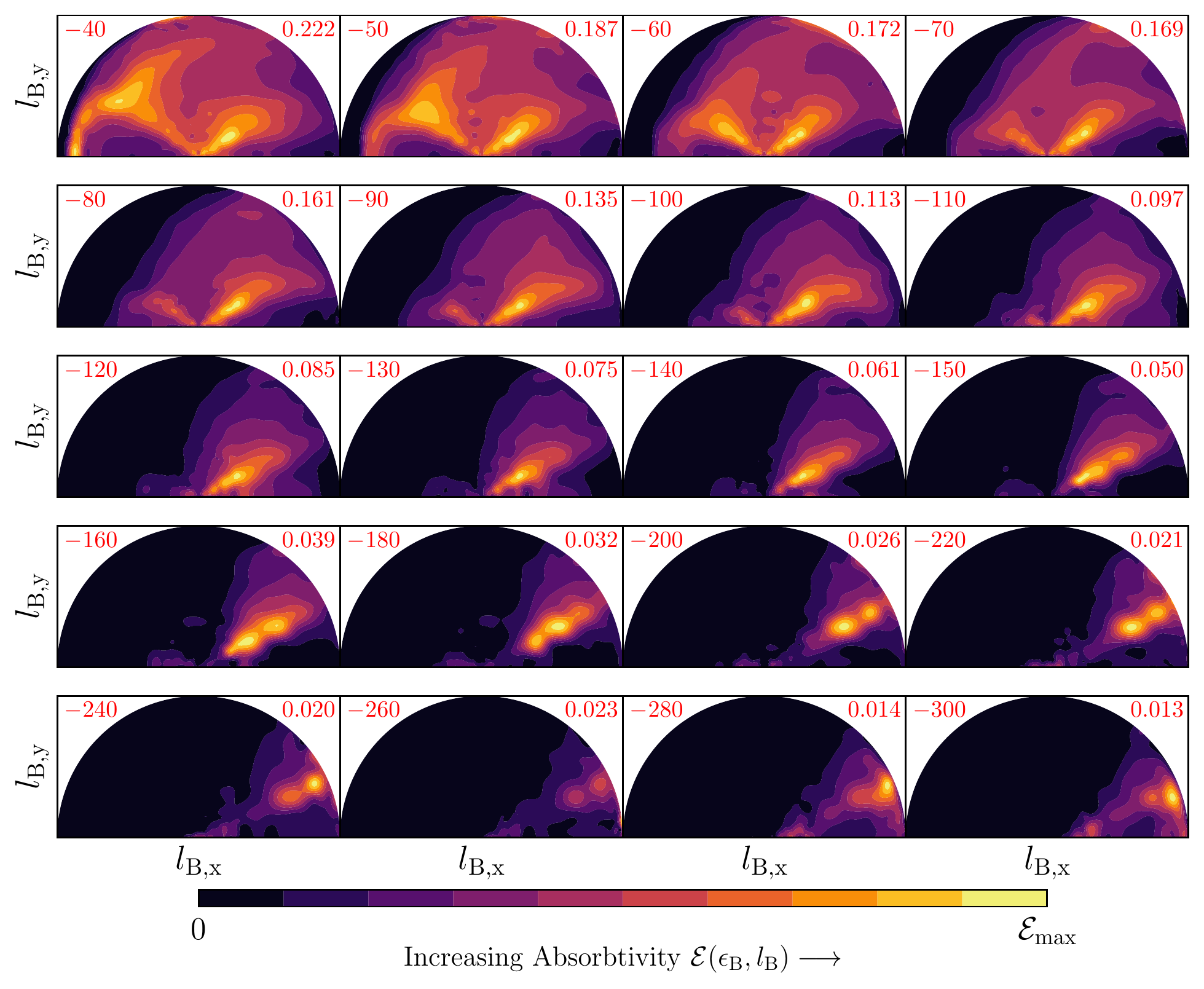}
    \caption{Chaotic absorptivity $\calE$ measured as function of $l_{\rm B,x},\, l_{\rm B,y}$ and $\epsilon_{\rm B}$. Note that $l_{\rm B,x},\, l_{\rm B,y}$ are the components of the binary angular momentum, and $\epsilon_{\rm B}$ is the binary energy. $\calE \left( l_{\rm B,x},\, l_{\rm B,y},\, \epsilon_{\rm B} \right)$ is presented as a sequence of colored contour plots, with each 2D slice corresponding to a fixed value of $\epsilon_{\rm B}$, shown on the top-left corner of each panel. To be able to resolve the absorptivity structures in each panel, the color scaling for each panel is scaled relative to its maximum absorptivity value, $\mathcal{E}_{\rm max}$, shown in top-right corner of each panel. To better see how these absorptivity structures evolve as a function of $\eps_{\rm B}$, we prepared a video out of the sequence of such slices. See YouTube or the following \href{http://phys.huji.ac.il/~barak_kol/resrch_supp/3body/3d_calE.mp4}{link}.$^1$} 
    \label{fig:measured_calE}
\end{figure*}
\footnotetext[1]{On YouTube search for ``Measured three-body chaotic absorptivity'', or download video from \url{http://phys.huji.ac.il/~barak_kol/resrch_supp/3body/3d_calE.mp4}.}

% We describe the simulations in the previous section. Furthermore, we have 

% \subsection{Direct measurement through three-body scattering experiments }
% \label{ssec:direct_cale}

\begin{figure}
    \centering
    \includegraphics[width = \columnwidth]{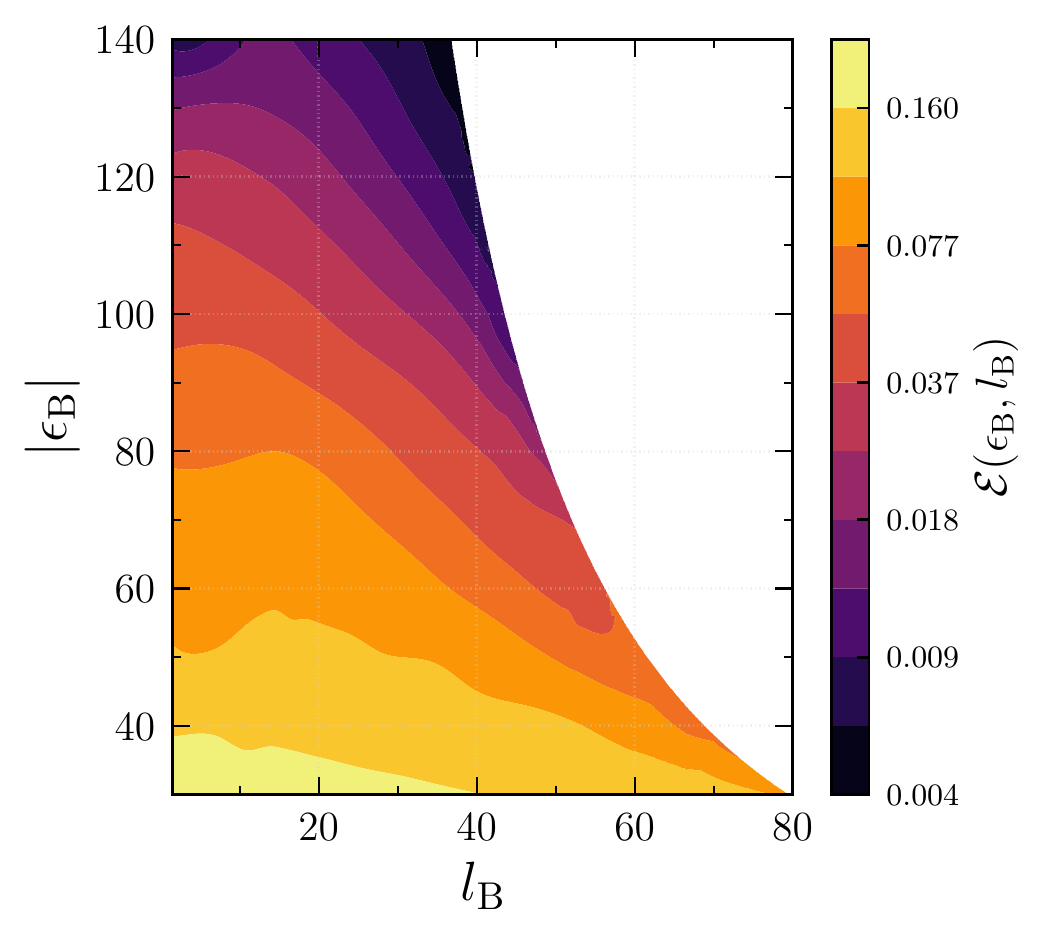}
    \caption{A contour plot showing the measured bi-variate chaotic absorptivity $\mathcal{E}$ as a function of binary energy $\eps_{\rm B}$ and binary angular momentum $l_{\rm B}$. Note that the contours are logarithmically spaced in $\calE$. }
    \label{fig:marginalized_calE}
\end{figure}

    % Marginalized measured chaotic absorptivity $\mathcal{E}$. The tri-variate function $\calE = \calE \left( l_{\rm B,x},\, l_{\rm B,y},\, \epsilon_{\rm B} \right)$ shown in Fig. \ref{fig:measured_calE}, is marginalized over the direction of $\vec{l}_B$ to produce $\calE = \calE \left( l_{\rm B},\, \epsilon_{\rm B} \right)$, where $l_{\rm B}$ is the magnitude of the binary angular momentum.
    %The contour plot for \textit{measured} 2d marginalized chaotic absorptivity $\mathcal{E}$ as a function of binary energy $|\epsilon_{\rm B}|$ and binary angular momentum $l_{\rm B}$.  Note that the distribution is marginalized as we average over different inclinations between $\vec{L}$ and $\vec{l_{\rm B}}$. 

\begin{figure*}
    \centering
    \includegraphics[width = \textwidth]{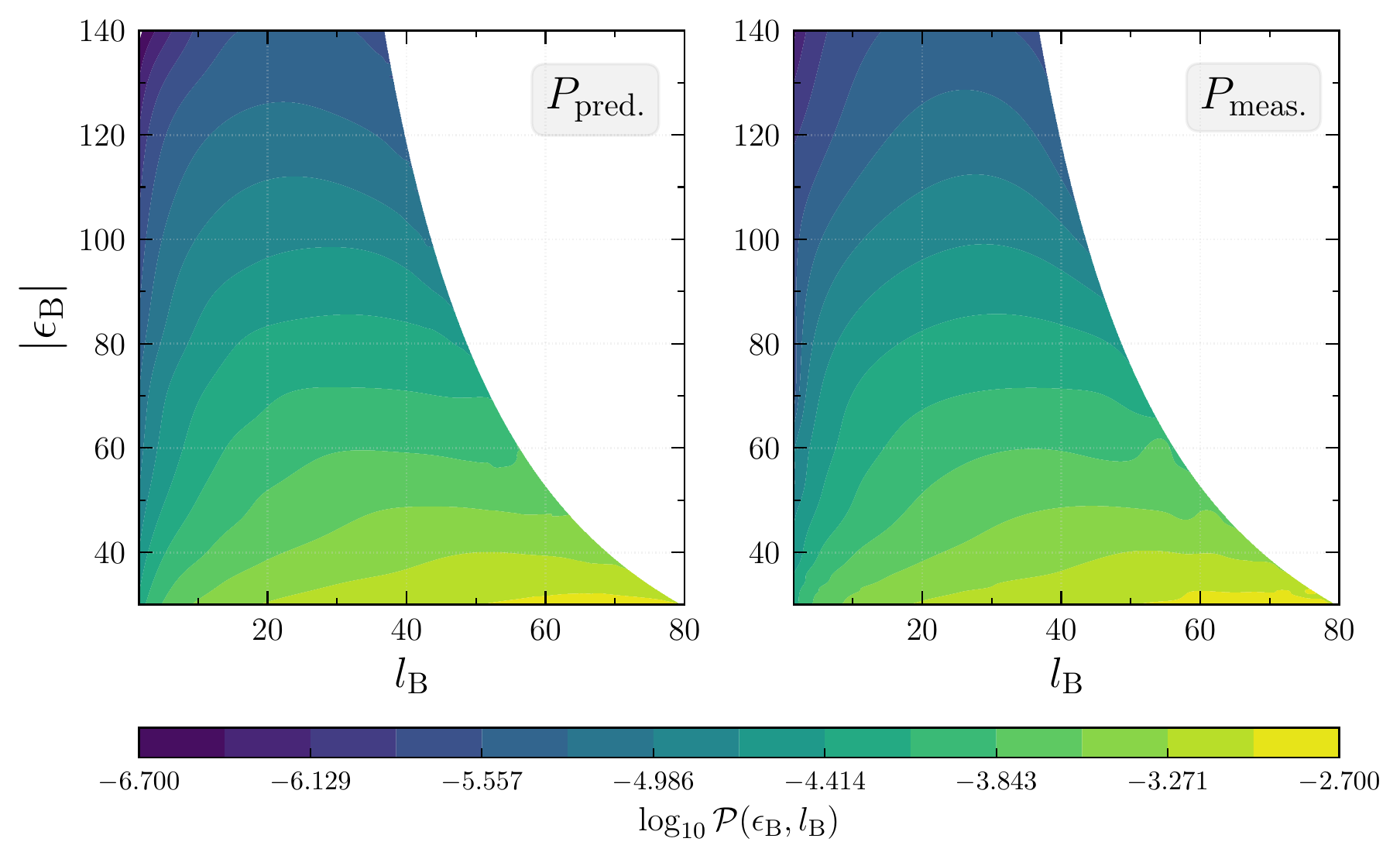}
    \caption{\textit{Predicted} (left) and \textit{measured} (right) chaotic outcome distributions as a function of binary energy $\epsilon_{\rm B}$ and binary angular momentum $l_{\rm B}$. The measured outcome distribution ($P_{\rm meas.}$) is obtained by running an ensemble of three-body interactions and filtering for chaotic escapes. The predicted outcome distribution ($P_{\rm pred.}$) is obtained by computing the chaotic absorptivity $\mathcal{E}$ (through our binary-single scattering experiments described in Section~\ref{sec:data}) and then multiplying it by the flux factor \eqref{dF} to get the above prediction. The predicted distribution has been scaled by a normalization factor so that the the median probability value in above region is same as measured distribution. Refer to Section~\ref{ssec:bivariate} for more details on this normalization. There is a striking resemblance between the measured and theoretically predicted chaotic outcome distribution. }
    \label{fig:predicted_vs_measured}
\end{figure*}

% Maybe a side by side comparison of the two outcome distributions? Maybe replace the 2d histogram with the other 2d distribution. What normalization to use? A fitted normalization is easier, because measuring the area of predicted absorptivity is not possible (we have limited phase space exploration also). The 2d rectangular grid is incompatible with the curved boundary and hence the cell has a wrong area and hence the issue. I can fix this by doing the correct normalization of area. Need to figure out what fraction of the cell is crossed by the boundary and then maybe can fix it?? Not sure...   

\section{Absorptivity Results}
\label{sec:data}

In this section, we discuss the direct measurements of chaotic absorptivity $\calE$ using simulations (as described in Section~\ref{sssec:sims_absorb}) and comparisons between the theoretically predicted and measured chaotic outcome distribution.

\subsection{Tri-variate absorption}
\label{ssec:trivariate}

The tri-variate absorptivity is the absorptivity $\calE$ as a function of three variables, namely, the binary energy $\eps_{\rm B}$, and the two components of the binary angular momentum along the total angular momentum vector ($l_{\rm B,x}$ and $l_{\rm B,y}$), where $l_{\rm B,x}$ is aligned with $\vec{L}$. A visualization of this three-dimensional space is provided in Figure~\ref{fig:3dsec}. 

To study the tri-variate absorption, we make various measurements of $\calE$ on a three dimensional grid as described in Section~\ref{sssec:sims_absorb}. We then interpolate on the grid across disks (a disk corresponds to a slice of constant $\eps_{\rm B}$) to compute a map of $\calE$ values for various binary energies $\epsilon_{\rm B}$. Figure~\ref{fig:measured_calE} shows these maps at various binary energies ranging from $\epsilon_{\rm B} = -40$ to $\epsilon_{\rm B} = -300$. The darker colors correspond to regions of low absorptivity, while brighter colors correspond to regions of high absorptivity. Note that the color range for each panel has been adjusted so that the brightest color corresponds to the maximum value of absorptivity ($\calE_{\rm max}$) in that panel. The color-bar has linear scaling in absorptivity in the range $[0,\calE_{\rm max}]$. For reference, the maximum value of absorptivity in each panel is denoted in the top right corner of respective panel. The binary energy value corresponding to the disk in each panel is denoted in the top left corner of respective panel. 

Looking at Figure~\ref{fig:measured_calE}, we observe a few patterns. The binary energy, not only influences the maximum value of $\calE$, but also the overall distribution of $\calE$ values on the disk.

Firstly, as we go from high binary energy values ($\epsilon_{\rm B} = -40$) to low binary energy values ($\epsilon_{\rm B} = -300$), a larger and larger fraction of the disk is \textit{not} absorbed into chaotic motion (that is, low absorptivity). This is because as the binary pair gets tighter (smaller $\epsilon_{\rm B}$), it becomes increasingly difficult for the incoming single particle to disrupt it into a state of rough equipartition, that is, into chaotic motion. Also, as the binary gets tighter, it carries less angular momentum $l_{\rm B}$. Thus, the single carries a larger fraction of the total angular momentum and hence will have a larger impact parameter, resulting in a higher probability of non-chaotic interaction. These non-chaotic interactions are either fly-bys or prompt exchanges. 

Secondly, as the binary gets tighter, the regions of zero absorptivity (dark color) increase in size from negative to positive $l_{\rm B,x}$ (from right to left). This is due to the retrograde (counter-rotating) and prograde (co-rotating) motion of the single relative to the motion of the binary. For $\eps_{\rm B} \gtrsim -120$, both prograde and retrograde encounters appear to have positive absorptivity, albeit only in the slightly misalligned case ($\sigma \approx 20^\circ$), because the coplanar case ($\sigma \approx 0^\circ$ or $180^\circ$) mostly consists of regular motion. For binaries with the smallest energies (e.g., $\eps_{\rm B} \lesssim -200$), the absorptivity is generally very small and is only non-zero in the case of prograde motion (positive $l_{\rm B,x}$). In general, the plot seems to show that chaotic scattering is more likely to occur in cases when motion of single is prograde relative to binary motion. 

% If motion of the single is retrograde (counter-rotating) relative to binary motion, then it is easier to have a regular scattering and hence low absorptivity. \AAT{(I don´t understand why it should be easier, or at least it is not explained. or are you just descriving the plot?)}

Thirdly, the overall structure of the $\calE$ space is quite complex as seen in Figure~\ref{fig:measured_calE}. For instance, one observes these two loci of high absorptivity in $\eps_{\rm B,x} = -40$ and these loci move through the $l_{\rm B,x} = l_{\rm B,y} = 0$ point with decreasing $\epsilon_B$. Note that the $l_{\rm B,x} = l_{\rm B,y} = 0$ point is the case where the binary has zero angular momentum and the single has all the angular momentum. It is beyond the scope of this paper to physically interpret the origin of these absorptivity structures and their evolution with $\eps_{\rm B}$.

% The main takeaway points of this discussion would be : i) to first order, the average absorptivity decreases with decreasing binary energy $\eps_{\rm B}$ and ii) there is a complicated dependence of absorptivity with the angular momentum of binary and single. 

\subsection{Bi-variate absorption and predicted outcome distribution}
\label{ssec:bivariate}

The bi-variate absorptivity is the absorptivity $\calE$ as a function of the binary energy $\eps_{\rm B}$ and binary angular momentum $l_{\rm B}$. As described in Figure~\ref{fig:pap_schema}, the motivation for measuring the bi-variate absorptivity is to be able to predict the outcome distribution for chaotic three-body interactions and thus compare directly to simulations. In this section, we describe our measurement of the bi-variate absorptivity and the resulting prediction of the outcome distribution.

The bi-variate absorptivity distribution can be visualized by condensing the concentric disks of equal $\eps_{\rm B}$ along lines of constant $l_{\rm B} = \sqrt{l_{\rm B,x}^2 + l_{\rm B,y}^2}$ in tri-variate absorptivity space (see Figure~\ref{fig:3dsec}). As described in Section~\ref{sssec:sims_absorb}, to compute bi-variate absorptivity maps, we run ensembles of simulations on a two-dimensional grid in $\eps_{\rm B}$ - $l_{\rm B}$ space. For a consistency check on our simulations, we confirmed that the bi-variate absorptivity measurements on this two dimensional grid are consistent with the bi-variate absorptivity measurement by marginalizing over the tri-variate absorptivity we measured in Section~\ref{ssec:trivariate}. After performing a 2D interpolation across this grid, we show the bi-variate absorptivity map in Figure~\ref{fig:marginalized_calE}. Note that the contours are logarithmically spaced, however, the corresponding color bar labels are shown in linear values for reading convenience. As observed in the tri-variate absorptivity map (see Section~\ref{ssec:trivariate}), the absorptivity decreases with decreasing binary energy $\eps_{\rm B}$ (or increasing $|\eps_{\rm B}|$). At a fixed $\eps_{\rm B}$, a dependency on the binary angular momentum magnitude is also observed, which is reflective of the complicated structures seen in Figure~\ref{fig:measured_calE}.

Given the absorptivity function $\calE(\eps_B,l_B)$, the flux-based statistical theory predicts the following outcome distribution \be
d \Gamma_s(\eps_B,l_B)_{\rm pred.} \propto \calE(\eps_B,l_B) \cdot \frac{dF}{d\eps_B\, dl_B} \; d\eps_B\, dl_B
\ee
 where the marginalized flux is given by \eqref{marg_flux}.
%\end{equation}
%In the flux-based statistical approach, with a measurement of the absorptivity $\calE$, we can compute the outcome distribution by \textcolor{red}{Viraj : I am not sure about the correct notation here to denote the outcome distribution. Could you please confirm this paragraph and also if it makes sense?}
%\begin{equation}
%    d \Gamma_s(u) \propto \calE(u) \cdot \frac{l_{\rm B}}{\eps_{\rm B}^{3/2}}
%\end{equation}
%where $u$ denotes the set of outcome parameters. The above factor is obtained by integrating the asymptotic flux $\mathcal{F}$ radially over $d^3 l_B$ (see Eqn~\ref{dF}).

Thus, we compute the predicted outcome distribution by: i) Taking the absorptivity grid in Figure~\ref{fig:marginalized_calE} and multiplying the absorptivity value at each grid point by the factor $l_{\rm B} / (-\eps_{\rm B})^{3/2}$ and then ii) performing a 2D spline approximation across this grid to obtain our prediction for the outcome distribution. Note that in multiplying the absorptivity distribution by this factor, the resulting outcome distribution will not be normalized and we will have to explicitly normalize it by a constant. However, in this scenario we cannot normalize this distribution the way we normalize any probability distribution by ensuring the total area under curve is one. That is because our $\eps_{\rm B}$ - $l_{\rm B}$ grid over which we compute the absorptivity (and hence the outcome distribution) is limited to only a certain part of the space : $\eps_{\rm B} \in [-150,-30]$ and $l_{\rm B} = [1.5,70]$. To enable a comparison with the measured outcome distribution, which is properly normalized, we scale our predicted outcome distribution such that its median value in the above region matches the median value of the measured outcome distribution in the same region.
Note that the precise value of normalization is not of importance here, but rather a comparison of the relative structures/probabilities in the outcome distribution. Following these steps, the left panel in Figure~\ref{fig:predicted_vs_measured} shows our predicted outcome distribution. 

% For comparison, the right panel in Figure~\ref{fig:predicted_vs_measured} shows the measured outcome distribution. In the following section, we will compare in detail these two quantities. 

\subsection{Comparison with measured outcome distribution}

As described in Section~\ref{sssec:outcome_meas}, we measure the chaotic outcome distribution by first running an ensemble of $10^7$ three-body interactions and then applying appropriate cuts to isolate the chaotic escapes. The right panel in Figure~\ref{fig:predicted_vs_measured} shows the contour plot for our measured outcome distribution. We obtain that plot by first binning the measured outcome distribution into a normalized 2D histogram (see Figure~\ref{fig:outcome_hist}). One notices that in Figure~\ref{fig:outcome_hist}, bins along the boundary of the histogram have anomalously low probabilities relative to their immediate neighbors not located on the boundary. As a reminder, the boundary occurs due to a forbidden region in $\epsilon_{\rm B} - l_{\rm B}$ space (see Equation~\ref{eqn:max_lb_eqn}). The reason for this boundary effect in the 2D histogram is that the rectangular binning of our histogram is not compatible with the shape of the boundary. Hence, when computing the probability density in each bin, the bin area goes outside the boundary, while the bin events/counts are purely located within the boundary, resulting in the anonymously low probability density relative to its non-boundary neighbors. We account for this effect by recomputing the density in each bin by considering only the bin area that is within the boundary. Note that this is not a perfect solution, however, it helps reduce the impact of the boundary. After accounting for this, we then smooth this histogram by 2D interpolation to get the contour plot shown in the right panel in Figure~\ref{fig:predicted_vs_measured}. 

\begin{figure*}
    \centering
    \includegraphics[width = \textwidth]{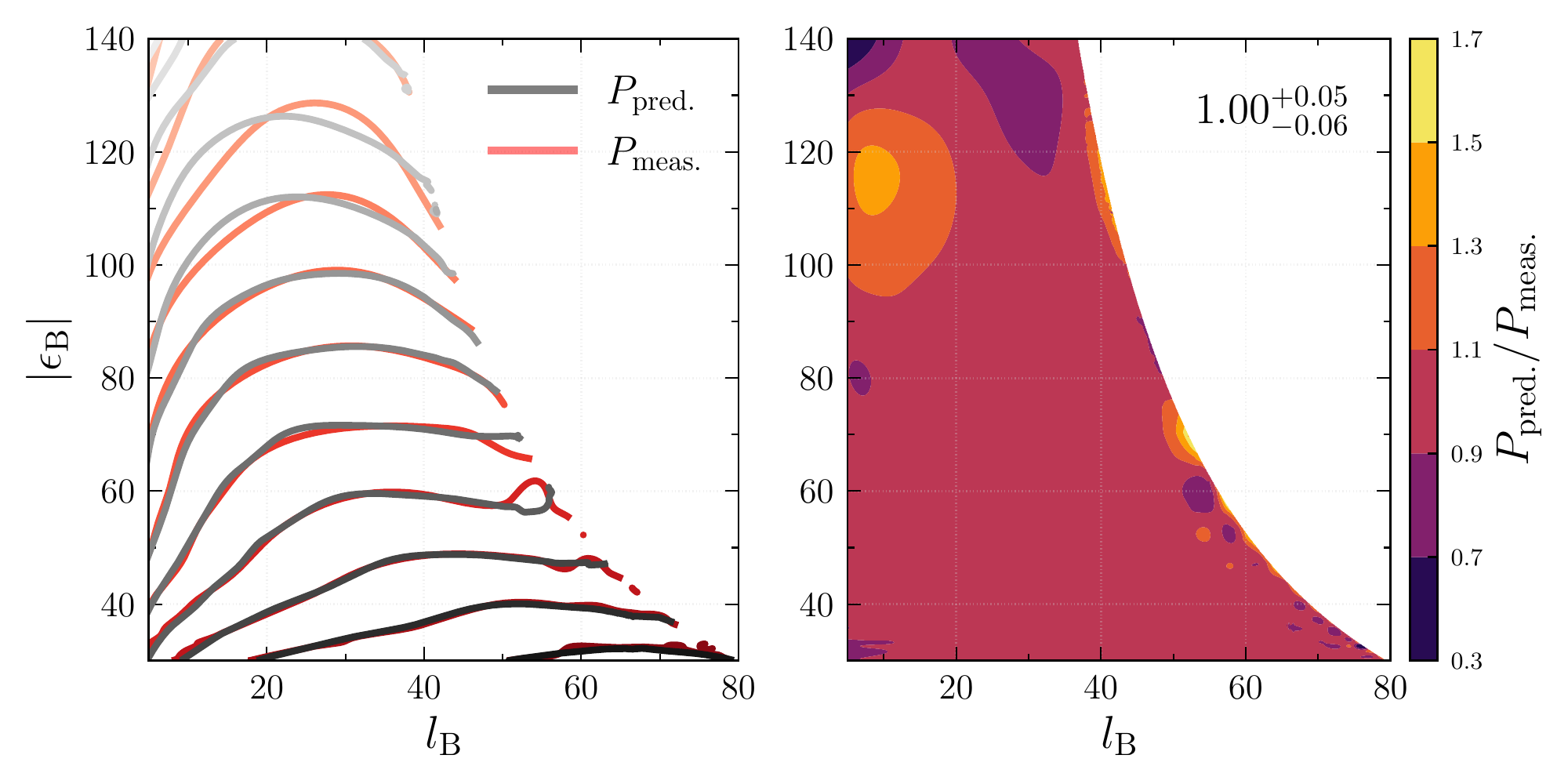}
    \caption{ Direct comparisons between the predicted ($P_{\rm pred.}$) and measured ($P_{\rm meas.}$) outcome distributions.  (\textit{Left Panel}) Overlaid contour plot of $P_{\rm pred.}$ and $P_{\rm meas.}$. Note that the contours are the same as the ones shown in Figure~\ref{fig:predicted_vs_measured}. (\textit{Right Panel}) Colored contour plot of the ratio $P_{\rm pred.}/ P_{\rm meas.}$. The top right corner shows the 16th and 84th percentile ranges of the ratio values in the shown region.  }
    \label{fig:residuals}
\end{figure*}

Figure~\ref{fig:predicted_vs_measured} shows a side-by-side comparison between our predicted outcome distribution (left panel) and the measured outcome distribution from simulations (right panel). Our predictions look quite similar to the measurement from simulations. For a direct comparison, the left panel in Figure~\ref{fig:residuals} shows the contours from these two distributions overlaid on top of each other. Overall good agreement is seen between the two, however with agreement deteriorating in regions with low statistics like high $|\eps_{\rm B}|$ and low $l_{\rm B}$. This is most likely due to our measured outcome distribution being impacted by a lack of robust interpolation in low statistics regions. Another way to compare the two is to look at the residual ratios of these two distributions. The right panel of Figure~\ref{fig:residuals} shows the contour plot for the ratio of our predicted outcome distribution to the measured outcome distribution ($P_{\rm pred.} / P_{\rm meas.}$). If the flux-based theory is valid, then throughout the $\eps_{\rm B} - l_{\rm B}$ space, this ratio should be close to $1$. As the plot shows, there is good agreement (${\sim} 10\%$ around 1) over most of the plane except for areas where we are impacted by low statistics in the $P_{\rm meas.}.$ The 16th and 84th percentile range of this contour plot is $[0.94,1.05]$ yielding a $\sim 6\%$ error estimate .

In conclusion, % within the limit of how well we can measure the outcome distribution without being impacted by low statistics, 
 in regions with good statistics, 
we find good agreement between the flux-based theory predictions for the chaotic outcome distribution and the corresponding measurement from simulations.

%%%%%%%%%%%%%%%%%%%%%%%%%%%%%%%%%%%%%%%%%%%%%%%%%%%%%%%%
\section{Summary and Discussion} 
\label{sec:summary}

% Results
% Performed a first measurement of the chaotic absorptivity function. Fixed equal masses. averaged over phases only and determined as a function the three conserved quantities ... . The working criterion for chaotic absorptivity is detailed in the corresponding subsection in section 3, as well as in section 4?
% Obtained the predicted outcome distribution using the flux-based theory (1,5).
% Improved the measurement of the chaotic outcome distribution performed in MKTL21. 
% Compared the predicted and measured chaotic outcome distributions and found them to be in very good agreement.

% Discussion and open Q
% Formulate a model for the absorptivity function. Ginat Perets.
% Extend the measurement of outcome statistics to be tri-variate and thereby match the measurement of calE. 
% Extend the measurement of calE: as a function of more variables or for unequal masses.
%% Clean the issue of the minimum 80 yr lifetime cut. At the moment it is not clear why this could be necessary.

\noindent {\bf Summary of results}. \begin{itemize}

\item We have performed a first measurement of $\calE$, the chaotic absorptivity function defined in \cite{Kol2020}. For concreteness, we fixed the masses to be equal, and measured $\calE$ as a function of those of its variables that are asymptotically conserved quantities, thereby averaging over the two pericenter angles. The values are presented in Fig. \ref{fig:measured_calE}. Our working criterion for chaotic absorption is detailed in Sect. \ref{sec:method}.  

\item We obtained the outcome distribution predicted by the flux-based theory through (\ref{dGamma},\ref{dF},\ref{marg_flux}).

\item We improved the measurement of the chaotic outcome distribution performed in MKTL21 in two ways : i) increased the number of three-body equal mass systems  from a million to 10 million realizations, and ii) adjusted the criterion for a chaotic escape, detailed in Sect. \ref{sec:data}, to match the criterion used in the direct measurement of chaotic absorptivity in this paper.
% (\textcolor{red}{Viraj: run on sentence and possible confusion with bracket so make it more clear.})
% BK: is this a correct description of the improvements? 

\item Finally, we compared the predicted and measured chaotic outcome distributions, and found them to be in good agreement, see figures \ref{fig:predicted_vs_measured},\ref{fig:residuals}.

\end{itemize}

\noindent {\bf Discussion}. The agreement of predicted and measured outcome distributions was found to 
% be within 10\% almost everywhere within 
 hold up to a $\sim 6\%$ error estimate throughout 
 the 2d parameter space. We consider this to be a detailed and precise agreement. This is more so given that this agreement is based on two rather different measured quantities: the absorptivity and the outcome distribution. We believe that the residuals arise from errors introduced by the smoothing procedure of the outcome distribution, especially near the boundary of parameter space and by lower outcome statistics at large $\left| \epsilon_{\rm B} \right|$ and low $l_{\rm B}$.

\presub {\bf Open questions}. \begin{itemize}

\item Formulate analytic models of the chaotic absorptivity function. Such models could partially replace the measurement that we performed. Inspection of Fig. \ref{fig:measured_calE} suggests that it would not be an easy task. On the other hand, the agreement with outcome distributions found in \cite{ginat20} should be interpreted in the current context as a proposal for such an analytic model of $\calE$. An analytic model of $\calE$ would open the way to improve analytical models of binary-single scattering, see \cite{2022MNRAS.517.3838L,2023MNRAS.519L..15G} for recent results on this topic.

\item Extend the measurement of outcome statistics to be tri-variate and thereby match the measurement of $\calE$. 

\item Extend the measurement of $\calE$. It could be measured as a function of one or two of the pericenter angles. In addition,  it could be extended to unequal masses.

\end{itemize}

\noindent {\bf In summary}, we have demonstrated the validity of the reduction of the three-body outcome statistics introduced by the flux-based statistical theory \cite{Kol2020}.

%%%%%%%%%%%%%%%%%%%%%%%%%%%%%%%%%%%%%%%%%%%%%%%%%%%%%%%%
\section*{Acknowledgements}

A.A.T. acknowledges support from JSPS KAKENHI Grant Number 21K13914 and from Horizon 2020 in the form of a Senior INTERACTIONS COFUND Fellowship. BK was partially supported by the Israel Science Foundation (grant No. 1345/21). Analyses presented in this paper were greatly aided by the following free software packages: \texttt{NumPy} (\citealt{numpy}), \texttt{Matplotlib} (\citealt{matplotlib}), \texttt{emcee} (\citealt{emcee}) and \texttt{Jupyter} (\citealt{jupyter}). This research has made extensive use of NASA\textquotesingle s Astrophysics Data System and arXiv. 
% Thank Nathan Leigh and Nicholas Stone for previous collaboration and discussions?

%%%%%%%%%%%%%%%%%%%%%%%%%%%%%%%%%%%%%%%%%%%%%%%%%%
\section*{Data Availability}

The \textsc{tsunami} code, the initial conditions and the simulation data underlying this article will be shared on reasonable request to the corresponding author. 

%%%%%%%%%%%%%%%%%%%% REFERENCES %%%%%%%%%%%%%%%%%%

% The best way to enter references is to use BibTeX:

\bibliographystyle{mnras}
\bibliography{main} % if your bibtex file is called example.bib

%%%%%%%%%%%%%%%%%%%%%%%%%%%%%%%%%%%%%%%%%%%%%%%%%%

%%%%%%%%%%%%%%%%% APPENDICES %%%%%%%%%%%%%%%%%%%%%

\appendix

% Don't change these lines
\bsp	% typesetting comment
\label{lastpage}
\end{document}